\documentclass[12pt,leqno]{article}
\usepackage{amsmath,amssymb}
\numberwithin{equation}{section}
\allowdisplaybreaks

\newtheorem{defin}{Definition}[section]
\newtheorem{prop}{Proposition}[section]
\newtheorem{corol}{Corollary}[section]

\newtheorem{rem}{Remark}[section]

\begin{document}
\title{Transitive Courant algebroids}
\author{{\small by}\vspace{2mm}\\Izu Vaisman}
\date{}
\maketitle
{\def\thefootnote{*}\footnotetext[1]%
{{\it 2000 Mathematics Subject Classification: 53D17}.
}}
\begin{center} \begin{minipage}{12cm}
A{\footnotesize BSTRACT. We express any Courant algebroid bracket
by means of a metric connection, and construct a Courant algebroid
structure on any orthogonal, Whitney sum $E\oplus C$ where $E$ is
a given Courant algebroid and $C$ is a flat, pseudo-Euclidean
vector bundle. Then, we establish the general expression of the
bracket of a transitive Courant algebroid, i.e., a Courant
algebroid with a surjective anchor, and describe a class of
transitive Courant algebroids which are Whitney sums of a Courant
subalgebroid with neutral metric and Courant-like bracket and a
pseudo-Euclidean vector bundle with a flat, metric connection. In
particular, this class contains all the transitive Courant
algebroids of minimal rank; for these, the flat term mentioned
above is zero. The results extend to regular Courant algebroids,
i.e., Courant algebroids with a constant rank anchor. The paper
ends with miscellaneous remarks and an appendix on Dirac linear
spaces.}
\end{minipage}
\end{center}
\vspace{5mm}
%\noindent
%begin{center} %\section %\end{center}
\section{The basics of Courant algebroids}
The framework of this note is the $C^\infty$-category. In the
literature, there are two notions of a {\it Courant algebroid},
which include a skew-symmetric and a non-skew-symmetric bracket,
respectively. These notions are the result of an effort to unify
the {\it Courant bracket} and the {\it Manin bracket}
\cite{Liu}.

We start with the definition of a non-skew-symmetric Courant
algebroid \cite{{Royt},{SW}}, with the simplifications indicated
in \cite{U}.
\begin{defin}\label{Calg} {\rm
A {\it Courant algebroid} is a pseudo-Euclidean vector bundle
$(E\rightarrow M,g)$ ($g$ is a symmetric, non degenerate inner
product on $E$) with an {\it anchor morphism} $\rho:E\rightarrow
TM$ and a general, $\mathbb{R}$-bilinear product $\star:\Gamma
E\times\Gamma E\rightarrow\Gamma E$ ($\Gamma$ denotes spaces of
global cross sections) such that $\forall e,e_1,e_2,e_3\in\Gamma
E$ one has:
$$\begin{array}{l}
1)\hspace{1cm}(\rho e)(g(e_1,e_2))=g(e\star e_1, e_2)+g(e_1,e\star
e_2),\vspace{2mm}\\ 2)\hspace{1cm}e\star e
=\partial(g(e,e)),\vspace{2mm}\\ 3)\hspace{1cm} e_1\star(e_2\star
e_3)=(e_1\star e_2)\star e_3 +e_2\star(e_1\star e_3),\end{array}$$
where $\partial=(1/2)\sharp_g\circ^t\rho:T^*M \rightarrow \Gamma
E$ and
\begin{equation}\label{defpartial}
\partial f=\partial(df)=\frac{1}{2}\sharp_g\,^t\rho(df),\hspace{1cm}
\forall f\in C^\infty(M)/\end{equation} Here, $t$ denotes
transposition and $\sharp$ is the musical isomorphism defined like
in Riemannian geometry.}
\end{defin}

By polarization, property 2) is equivalent with the prescription
of the symmetric part of the product
\begin{equation}\label{simetrizare} (e_1,e_2)_{*}=
\frac{1}{2}(e_1\star e_2+e_2\star e_1)=
\partial(g(e_1,e_2)).\end{equation}
On the other hand, the definition of the operator $\partial$ is
equivalent with
\begin{equation}\label{partialprinrho} g(\partial f,e)=
\frac{1}{2}(\rho e)f,\hspace{5mm}\forall f\in C^\infty(M),\forall
e\in\Gamma E.\end{equation} Hence, property 2) is also equivalent
with \begin{equation}\label{partialsisim}
g(e,(e_1,e_2)_{*})=\frac{1}{2}(\rho e)(g(e_1,e_2)).\end{equation}
Thus, modulo property 1) we may replace 2) by a condition that
does not contain $\partial$ namely,
\begin{equation}\label{3farapartial}
2g(e,(e_1,e_2)_{*})= g(e\star e_1,e_2) + g(e_1,e\star e_2).
\end{equation}

Let us also consider the skew-symmetric part of the product:
\begin{equation}\label{antisimetrizare} [e_1,e_2] = [e_1,e_2]_{*}=
\frac{1}{2}(e_1\star e_2-e_2\star e_1).\end{equation} With
(\ref{simetrizare}) and (\ref{antisimetrizare}) we get
\begin{equation}
\label{generalC}e_1\star e_2=[e_1,e_2]+\partial
g(e_1,e_2).\end{equation}

From the properties postulated by Definition \ref{Calg} one can
deduce
\begin{prop}\label{Uchino} {\rm\cite{U}} Let $(E\rightarrow M,g,\rho,\star)$ be
a Courant algebroid. Then, $\forall e,e_1,e_2\in\Gamma E,\forall
f\in C^\infty(M)$, the following properties hold: $$\label{proU}
\begin{array}{l} a)\hspace{1cm} e_1\star (fe_2)=f(e_1\star e_2)+
((\rho e_1)f)e_2,\vspace{2mm}\\ b)\hspace{1cm} (fe_1)\star
e_2=f(e_1\star e_2)- ((\rho e_2)f)e_1 +2g(e_1,e_2)\partial
f,\vspace{2mm}\\ c)\hspace{1cm} (\partial f)\star e=0,\;\;
e\star(\partial f)=\partial((\rho e)f),\vspace{2mm}\\ d)
\hspace{1cm} \rho(e_1\star e_2)=[\rho e_1,\rho
e_2]_{TM},\vspace{2mm}\\ e) \hspace{1cm} \rho(\partial
f)=0.\end{array}$$ \end{prop} \noindent {\bf Proof.} Property a)
follows from the comparison of the results of expressing $(\rho
e)(g(fe_1,e_2))$ in two ways, first by applying property 1) of
Definition \ref{Calg} straightforwardly, second by using the
Leibniz rule for the vector field $\rho e$ applied to the product
$fg(e_1,e_2)$ and then property 1) for $g(e_1,e_2)$.

Property b) follows from a) by using (\ref{simetrizare}) and the
fact that, on functions, $\partial$ satisfies the Leibniz rule.

Notice that properties a) and b) show that a Courant algebroid
product is an operator of the local type (i.e., $(e_1\star
e_2)(x)$ depends only on the restrictions of $e_1,e_2$ to a
neighborhood of the point $x\in M$).

Now, let us denote \begin{equation}\label{obstructiaL}
\mathcal{L}(e_1,e_2,e_3)=e_1\star(e_2\star e_3)-
(e_1\star e_2)\star e_3-e_2\star(e_1\star e_3).\end{equation} From
(\ref{simetrizare}), it follows that
\begin{equation}\label{siminL} \mathcal{L}(e_1,e_2,e)
+ \mathcal{L}(e_2,e_1,e)=-2[\partial(g(e_1,e_2))\star e].
\end{equation} Since any function $f\in C^\infty(M)$ may locally
be written as $$f=g(e_1,\frac{f}{g(e_1,e_1)}e_1),$$ where $e_1$ is
not $g$-isotropic, equation (\ref{siminL}) and property 3) of
Definition \ref{Calg} imply the first part of c). The second part
of c) follows from the first part and the relation
(\ref{simetrizare}).

In order to get property d) (where the right hand side is a Lie
bracket of vector fields on $M$) we start with
(\ref{partialprinrho}), written as
$$(\rho e_2)f=2g(\partial f,e_2),$$ and apply $\rho e_1$, while using
1) of Definition \ref{Calg} and (\ref{partialprinrho}) again. The
result is
\begin{equation}\label{intermedresult}
(\rho e_1)(\rho e_2)f = 2g(e_1\star(\partial f),e_2)
+\rho(e_1\star e_2)f. \end{equation}In view of the second part of
c), this relation becomes
$$(\rho e_1)(\rho e_2)f=\rho(e_1\star e_2)f+ (\rho e_2)(\rho
e_1)f.$$ Since $f$ is an arbitrary function, we got precisely d).

Finally, from b) we get
$$\rho((fe_1)\star e_2)=f\rho(e_1\star e_2) - ((\rho e_2)f))\rho
e_1 + 2g(e_1,e_2)\rho(\partial f),$$ and, if we use d) in the two
sides of the previous relation while assuming $g(e_1,e_2)\neq0$,
we deduce the property e). Q.e.d.
\begin{rem}\label{obsdupaU} {\rm Property c) implies that the
skew-symmetric part of a Courant algebroid product satisfies the
property
\begin{equation}\label{cptantisim} [e,\partial f]=\frac{1}{2}
\partial((\rho e)f).\end{equation} Property e) is equivalent with
$g(\partial f,\partial g)=0$, $\forall f,g\in C^\infty(M)$ i.e.,
with the fact that $im\,\partial_x$ is a $g$-isotropic subspace of
the fiber $E_x$, $\forall x\in M$. We also note that property e)
is implied by b) and d). Finally, the computation used in the
proof of d) is reversible in the sense that d) and formula
(\ref{intermedresult}) imply property c) as well as its
consequence (\ref{cptantisim}).}
\end{rem}
\begin{rem}\label{cazulrhozero} {\rm
If $\rho=0$ we have $\partial=0$ and the Courant algebroid is just
a bundle of Lie algebras with a pseudo-Euclidean metric $g$ that
is invariant for the Lie algebra structure of each
fiber.}\end{rem}

From Definition \ref{Calg} and Proposition \ref{proU} it follows
that the skew-symmetric part (\ref{antisimetrizare}) of a product
$\star$ has the properties indicated by the following proposition.
\begin{prop} \label{propLiu} {\rm\cite{Royt}} For any Courant
algebroid and $\forall e,e_1,e_2,e_3\in\Gamma E,\,\forall f\in
C^\infty(M)$ one has:
$$\begin{array}{ll}
i)\hspace*{8mm}&\rho[e_1,e_2]=[\rho e_1,\rho e_2],\vspace{2mm}\\
ii)\hspace*{8mm}& im(\sharp_g\circ\,^t\hspace{-1mm}\rho)\subseteq
ker\,\rho\vspace{2mm}\\ iii)\hspace*{8mm}&
\sum_{Cycl}[[e_1,e_2],e_3]=
\frac{1}{3}\partial[\sum_{Cycl}g([e_1,e_2],e_3)],\vspace{2mm}\\
iv)\hspace*{8mm}& [e_1,fe_2]=f[e_1,e_2]+((\rho e_1)f)e_2 -
g(e_1,e_2)\partial f,\vspace{2mm}\\ v)\hspace*{8mm}& (\rho
e)[g(e_1,e_2)]=g([e,e_1]+\, \partial g(e,e_1),e_2) +
g(e_1,[e,e_2]+\, \partial g(e,e_2)).
\end{array}$$ \end{prop} {\bf Proof.} Except for iii), these properties
are immediate consequences of Definition \ref{Calg} and
Proposition \ref{proU}. For iii), let us denote
\begin{equation}\label{obstructiaJ}
\mathcal{J}(e_1,e_2,e_3) = [[e_1,e_2],e_3] + [[e_2,e_3],e_1] +
[[e_3,e_1],e_2]. \end{equation} Using (\ref{generalC}) we get
\begin{equation}\label{ajutor2} \mathcal{L}(e_1,e_2,e_3)=
\sum_{Cycl(1,2,-3)}\partial(g(e_1,[e_2,e_3])) \end{equation} $$+
\sum_{Cycl(1,-2,3)}[e_1,\partial(g(e_2,e_3))] -\frac{1}{2}
\sum_{Cycl(-1,2,3)}\partial((\rho e_1)(g(e_2,e_3))) -
\mathcal{J}(e_1,e_2,e_3),$$ where, in the sums of the right hand
side, the indices $(1,2,3)$ move cyclically while the signs of the
terms are as indicated in the summation index.

Now, taking the sum of (\ref{ajutor2}) over cyclic permutations of
$(1,2,3)$ and using (\ref{siminL}), we get
\begin{equation} \label{ajutor3}
6\mathcal{J}(e_1,e_2,e_3)=\sum_{Cycl(1,2,3)}\{\mathcal{L}(e_1,e_2,e_3)
- \mathcal{L}(e_2,e_1,e_3)\end{equation} \nopagebreak
$$+2\partial(g(e_1,[e_2,e_3]))\}.$$
Formula (\ref{ajutor3}) shows that property 3) of Definition
\ref{Calg} implies the present property iii). Q.e.d.

Following \cite{Liu}, one has \begin{defin} \label{skewCalg} {\rm
A {\it skew-symmetric Courant algebroid} is a pseudo-Eucli\-dean
vector bundle $(E\rightarrow M,g)$  with an anchor morphism
$\rho:E\rightarrow TM$, and a skew-symmetric bracket $[\;,\;]$ on
$\Gamma E$ such that properties i)-v) of Proposition \ref{propLiu}
hold.}\end{defin}
\begin{rem} \label{Uptantisim} {\rm \cite{U}} {\rm With the same proofs as for
Proposition \ref{proU}, we can see that property v) of Proposition
\ref{propLiu} implies iv) and properties i), iv) imply ii).
Therefore, conditions i), iii), v) suffice in the definition of a
skew-symmetric Courant algebroid.}\end{rem}
\begin{rem}\label{imrhofoliatie}  {\rm Just like in the case of a
Lie algebroid, for a skew-symmetric Courant algebroid
$\mathcal{D}=im\,\rho$ is a generalized foliation on $M$. Indeed,
from property i), Proposition \ref{propLiu}, it follows that
$\mathcal{D}=span\{\rho e\,/\,\forall e\in\Gamma E\}$ is spanned
by a Lie algebra, and, for any vector field $X=\rho e$, $\forall
x\in M$,
$dim\,\mathcal{D}_{(exp(tX))(x)}=rank\,\rho|_{(exp(tX))(x)}
=const.$ (For the latter assertion, look at a {\it Lie derivative}
of $\rho$ defined as for a tensor field and notice that this Lie
derivative vanishes.) Then the result follows from the
Sussmann-Stefan-Frobenius theorem (e.g., \cite{V94}, Theorem
2.9$''$).}\end{rem}

For simplicity, in what follows we will use the language provided
by the next definitions
\begin{defin} \label{ancoracourant} {\rm A {\it Courant anchor} is
a vector bundle morphism $\rho:E\rightarrow TM$, where $(E,g)$ is
a pseudo-Euclidean vector bundle over the manifold $M$, which is
such that $im\,\rho$ is a (generalized) foliation and the
corresponding morphism $\partial:T^*M\rightarrow E$ defined by
(\ref{defpartial}) has a $g$-isotropic image. A triple
$(E,g,\rho)$ where $\rho$ is a Courant anchor will be called a
{\it Courant vector bundle}. A Courant vector bundle endowed with
a skew-symmetric bracket $[\,,\,]$ on $\Gamma E$, which satisfies
properties v) and i) (therefore, also, iv) and ii) of a
skew-symmetric Courant algebroid will be called a {\it pre-Courant
algebroid}.}
\end{defin}

Notice that, since, $\forall x\in M$, both $im\,\rho$ and
$im\,\partial$ have the dimension equal to $rank\,\rho$, a Courant
vector bundle must satisfy the condition $rank\,\rho\leq
b\leq(1/2)rank\,E$, where $b$ is the smallest between the
positive-negative inertia indices of $g$. Furthermore, in view of
Remark \ref{obsdupaU}, formula (\ref{cptantisim}) holds for any
pre-Courant algebroid. Using (\ref{cptantisim}) it is easy to
check that any pre-Courant algebroid satisfies property iii) of
Proposition \ref{propLiu} if at least one of the arguments $e_a\in
im\,\partial$ $(a=1,2,3)$ (the cross sections of the subbundle
$im\,\partial$ are locally spanned over $C^\infty(M)$ by cross
sections of the form $\partial f$, $(f\in C^\infty(M))$).

By Proposition \ref{propLiu}, the skew-symmetric part of a Courant
algebroid product yields a skew-symmetric Courant algebroid
bracket. The converse is also true: \begin{prop}\label{propRoyt}
{\rm(E.g., \cite{Royt}).} If the bracket $[\,,\,]$ satisfies the
properties of a skew-symmetric Courant algebroid, the product
$\star$ defined by (\ref{generalC}) satisfies the properties of a
Courant algebroid.
\end{prop} \noindent {\bf Proof.} Obviously, 1) and 2) of Definition
\ref{Calg} hold. Moreover, as explained in Remark \ref{obsdupaU},
property c) of Proposition \ref{proU} holds independently of
formula (\ref{siminL}). Accordingly, now, (\ref{siminL}) proves
that $$\mathcal{L}(e_1,e_2,e_3) + \mathcal{L}(e_2,e_1,e_3)=0.$$ On
the other hand, if iii) of Proposition \ref{propLiu} holds,
formula (\ref{ajutor3}) implies
$$\mathcal{L}(e_1,e_2,e_3) - \mathcal{L}(e_2,e_1,e_3)=0.$$
Hence, we are done. Q.e.d.

The following proposition shows the possible changes of the
bracket of a pre-Courant algebroid.
\begin{prop}\label{strsugerate} Let $(E,g,\rho,[\,,\,])$ be a
pre-Courant algebroid. Then, the formula
\begin{equation}\label{sugerat1}
[e_1,e_2]'=[e_1,e_2]+\lambda(e_1,e_2),\hspace{5mm}e_1,e_2\in\Gamma
E,\end{equation} where $\lambda\in\Gamma(\wedge^2E^*\otimes E)$,
$\rho\circ\lambda=0$ and
\begin{equation}\label{Lambda}
\Lambda(e_1,e_2,e_3)=g(\lambda(e_1,e_2),e_3) \end{equation}
is totally skew-symmetric (i.e., $\Lambda\in\Gamma\wedge^3E^*$),
yields all the pre-Courant algebroid structures on $(E,g,\rho)$.
Furthermore, the bracket {\rm(\ref{sugerat1})} is that of a
Courant algebroid iff
\begin{equation}\label{condcociclu0}
(\partial_{[\,]}\Lambda)(e_1,e_2,e_3)=\mathcal{J}(e_1,e_2,e_3)
-\frac{1}{3}\partial[\sum_{Cycl}g([e_1,e_2],e_3)],\end{equation}
where
\begin{equation}\label{lambda[]}
(\partial_{[\,]}\Lambda)(e_1,e_2,e_3)=\partial(\Lambda(e_1,e_2,e_3))
-\sum_{Cycl(1,2,3)}\{ \lambda(\lambda(e_1,e_2),e_3)\end{equation}
$$+
\lambda([e_1,e_2],e_3)+[\lambda(e_1,e_2),e_3]\}\in
\Gamma(\wedge^3E^*\otimes E).$$
\end{prop}
\noindent {\bf Proof.} The difference of two brackets of
pre-Courant algebroid structures is a form
$\lambda\in\Gamma(\wedge^2E^*\otimes E)$ because of property iv),
Proposition \ref{propLiu}. The indicated conditions for $\lambda$
are equivalent with properties i) and v) of the same Proposition,
respectively. Notice that $\rho\circ\lambda=0$ is equivalent with
\begin{equation}\label{condiptLambda}
\Lambda(e_1,e_2,\partial f)=0,\hspace{1cm}\forall f\in
C^\infty(M),
\end{equation} therefore, in view of the skew-symmetry of
$\Lambda$, also equivalent with
\begin{equation}\label{condilambda} \lambda(e,\partial f)=0,
\hspace{5mm}\forall e\in\Gamma E,\;\forall f\in C^\infty(M).
\end{equation} Finally, a technical
computation shows that property iii) is equivalent with
(\ref{condcociclu0}), (\ref{lambda[]}). Q.e.d.
\begin{rem}\label{subalgisotropa} {\rm \cite{{C},{Liu}}
Any $g$-isotropic subalgebroid (i.e., a vector subbundle that is
closed by brackets) of a skew-symmetric Courant algebroid is a Lie
algebroid.}\end{rem}
\section{Courant brackets and metric connections}
We will get further insight into the structure of the bracket of a
Courant algebroid by using a metric connection $\nabla$ on the
pseudo-Euclidean bundle $(E,g)$, which means that
\begin{equation} \label{nablametric} X(g(e_1,e_2)) =
g(\nabla_Xe_1,e_2) + g(e_1,\nabla_Xe_2)\hspace{5mm}(X\in\Gamma
TM).\end{equation} If we also have the morphism $\rho:E\rightarrow
TM$, we define the {\it $\rho$-torsion} $T_{(\nabla,\rho)}
\in\Gamma(\wedge^2E^*\otimes TM)$ by the formula
\begin{equation}\label{pstorsiunepeE} T_{(\nabla,\rho)}(e_1,e_2)=
\rho(\nabla_{\rho e_1}e_2-\nabla_{\rho e_2}e_1)- [\rho e_1,\rho
e_2]. \end{equation}

Now, we can prove
\begin{prop}\label{crosetcuA} 1) The formula
\begin{equation}\label{solparticcugamma}
[e_1,e_2]_0=\nabla_{\rho e_1}e_2 -\nabla_{\rho e_2}e_1
-\gamma(e_1,e_2),
\end{equation} where $\gamma$ is defined by the equality
\begin{equation}\label{defluigamma}
g(\gamma(e_1,e_2),e)= \frac{1}{2}[g(e_1,\nabla_{\rho e}e_2) -
g(e_2,\nabla_{\rho e}e_1)],\hspace{2mm}\forall e\in\Gamma E,
\end{equation} defines a skew-symmetric bracket on $\Gamma E$ that
satisfies property v) of Proposition {\rm\ref{propLiu}}.

2) The most general bracket that satisfies v) is given by
\begin{equation}\label{solgencubeta}
[e_1,e_2]=[e_1,e_2]_0 - \beta(e_1,e_2),
\end{equation} where $\beta\in\Gamma(\wedge^2E^*\otimes E)$ and
\begin{equation} \label{definluiB} B(e_1,e_2,e_3)=
g(\beta(e_1,e_2),e_3) \end{equation} is totally skew-symmetric,
i.e., $B\in\Gamma\wedge^3E^*$.

3) The bracket {\rm(\ref{solgencubeta})} also satisfies property
i) of Proposition {\rm\ref{propLiu}} iff the following two
equalities hold \begin{equation}\label{rhopartialzero}
\rho\circ\partial=0,\end{equation}
\begin{equation}\label{condicubeta}
\rho(\beta(e_1,e_2))=T_{(\nabla,\rho)}(e_1,e_2),\hspace{1cm}\forall
e_1,e_2\in\Gamma E,\end{equation} where $T_{(\nabla,\rho)}$ is the
$\rho$-torsion of $\nabla$. In particular, if $(E,g,\rho)$ has a
metric connection $\nabla$ with zero $\rho$-torsion, the bracket
{\rm(\ref{solparticcugamma})} satisfies property i).\end{prop}
\noindent {\bf Proof.} First, let us notice that formula
(\ref{defluigamma}) actually defines $\gamma$ because if this
formula holds for $e\in\Gamma E$ it also holds for $fe$, $\forall
f\in C^\infty(M)$, and since $g$ is non degenerate. Then, a
technical calculation, which takes into account the metric
property (\ref{nablametric}) shows that the bracket
(\ref{solparticcugamma}) satisfies property v).

Now, the difference $\beta$ between two brackets that satisfy v)
is such that the corresponding $B$ given by (\ref{definluiB}) is
totally skew-symmetric, and $\beta\in\Gamma(\wedge^2E^*\otimes
E)$, respectively $B\in\wedge^3E^*$, because v) implies property
iv) of Proposition
\ref{propLiu}. These remarks justify the general formula
(\ref{solgencubeta}).

Finally, from (\ref{solgencubeta}), we get
\begin{equation}\label{ajutorpti}
\rho[e_1,e_2]= [\rho e_1,\rho e_2] + T_{\nabla,\rho}(e_1,e_2) -
\rho(\gamma(e_1,e_2)) - \rho(\beta(e_1,e_2)).
\end{equation} We know that, if properties i), v) of Proposition
\ref{propLiu} hold so does iv) and then ii), and ii) is equivalent
with (\ref{rhopartialzero}). Then, by using (\ref{defluigamma})
for $e=\partial f$, $f\in C^\infty(M)$, we get
$\rho\circ\gamma=0$, and (\ref{ajutorpti}) implies
(\ref{condicubeta}). The converse also is clear from
(\ref{ajutorpti}). The last assertion follows since, if
$T_{(\nabla,\rho)}=0$, we may chose $\beta=0$. Q.e.d.
\begin{rem}\label{condicuB} {\rm By applying the vector fields of
formula (\ref{condicubeta}) to an arbitrary function $f$, it
follows that (\ref{condicubeta}) is equivalent with
\begin{equation}\label{condicuB2} B(e_1,e_2,\partial
f)=\frac{1}{2}(T_{(\nabla,\rho)}(e_1,e_2))f,\hspace{5mm}e_1,e_2\in\Gamma
E,\;f\in C^\infty(M). \end{equation}}\end{rem}

Formula (\ref{solgencubeta}), where $B$ defined by
(\ref{definluiB}) is skew-symmetric and (\ref{rhopartialzero}),
(\ref{condicubeta}) hold, defines all the pre-Courant algebroid
brackets on a given Courant vector bundle. Among them, the Courant
algebroid brackets are obtained if condition iii) of Proposition
\ref{propLiu} also holds. In order to express the latter, let us
denote
\begin{equation}\label{defCrond} \mathcal{C}(e_1,e_2,e_3)=
\mathcal{J}(e_1,e_2,e_3)-\frac{1}{3}\partial[\sum_{Cycl(1,2,3)}
g([e_1,e_2],e_3)], \end{equation} where $\mathcal{J}$ is defined
by (\ref{obstructiaJ}), and denote by an index $0$ the same
expression for the bracket $[\,,\,]_0$. Assuming that property i)
holds, it follows that
\begin{equation}\label{Czero} \mathcal{C}_0(e_1,e_2,e_3)=
\sum_{Cycl(1,2,3)}\{(\nabla_{\rho e_3}\gamma)(e_1,e_2)+
\gamma(\gamma(e_1,e_2),e_3)\end{equation}
$$- R_\nabla(\rho e_1,\rho e_2)e_3\},$$

\begin{equation}\label{CprinCzero} \mathcal{C}(e_1,e_2,e_3)=
\mathcal{C}_0(e_1,e_2,e_3) +
\sum_{Cycl(1,2,3)}\{\beta(\beta(e_1,e_2),e_3)\end{equation} $$-
\beta([e_1,e_2]_0,e_3) -[\beta(e_1,e_2),e_3]_0\},$$ where
$R_\nabla$ is the curvature of the connection $\nabla$ and
$\nabla\gamma$ is defined as if $\gamma$ would be a tensor. These
formulas express the condition for the bracket
(\ref{solgencubeta}) to satisfy property iii), $\mathcal{C}=0$, by
means of the metric connection $\nabla$. In particular, if there
exists a metric connection of zero $\rho$-torsion, the
corresponding bracket (\ref{solparticcugamma}), which satisfies
i), also satisfies iii) iff $\mathcal{C}_0=0$, and this condition
reminds us the Bianchi identity for a linear connection with
torsion on a differentiable manifold $M$.

As an application of the results given in this section we have
\begin{prop}\label{strsugerate0} Let $(E,g,\rho,[\,,\,])$ be a
pre-Courant algebroid. Then, for any pseudo-Euclidean vector
bundle $(C,g_0)$ over $M$ and any metric connection $\nabla$ on
$C$, the brackets
\begin{equation}\label{sugerat2}
[e_1,e_2],\;[e,c]=-[c,e]=\nabla_{\rho
e}c,\;[c_1,c_2]=-\gamma_0(c_1,c_2),\end{equation} where
$e,e_1,e_2\in\Gamma E$, $c,c_1,c_2\in\Gamma C$ and
$\gamma_0(c_1,c_2)\in\Gamma E$ is defined by
\begin{equation} \label{gammazero}
g(\gamma_0(c_1,c_2),e)=\frac{1}{2}[g_0(c_1,\nabla_{\rho e}c_2) -
g_0(c_2,\nabla_{\rho e}c_1)], \end{equation} define a pre-Courant
algebroid structure on $(E\oplus C,g\oplus g_0,\rho\oplus 0)$. If
the original algebroid is a Courant algebroid,
{\rm(\ref{sugerat2})} yields a Courant algebroid structure iff the
connection $\nabla$ is flat.\end{prop} \noindent {\bf Proof.}
Straightforward checks show that the brackets (\ref{sugerat2})
satisfy property v) of a pre-Courant algebroid. For a full
justification of property i) we must also notice that
$\rho(\gamma_0(c_1,c_2)=0$. This follows since (\ref{gammazero})
implies $g(\gamma_0(c_1,c_2),\partial f)=0$, $\forall f\in
C^\infty(M)$. For the last assertion of the proposition we refer
to property iii) of a Courant algebroid. For arguments
$e_1,e_2,e_3$, iii) holds if $E$ is a Courant algebroid and for
arguments $c_1,c_2,c_3$ iii) follows from
$\rho(\gamma_0(c_1,c_2))=0$. For arguments $e_1,e_2,c$, iii) is
just $R_\nabla(\rho e_1,\rho e_2)c=0$, where $R_\nabla$ is the
curvature of $\nabla$. Finally, for arguments $e,c_1,c_2$, iii)
means \begin{equation}\label{cazulecc} [e,\gamma_0(c_1,c_2)] -
\gamma_0(\nabla_{\rho e}c_1,c_2) - \gamma_0(c_1,\nabla_{\rho
e}c_2)\end{equation}
$$-\frac{1}{2}\partial\{g_0(\nabla_{\rho e}c_1,c_2) -
g_0(c_1,\nabla_{\rho e}c_2)\}=0.$$ Since for a pre-Courant
algebroid property iv) holds, the left hand side of the previous
equality is $C^\infty(M)$-linear and it suffices to check it for a
local basis of $C$. If connection $\nabla$ is flat, $C$ has local
$\nabla$-parallel bases and (\ref{cazulecc}) obviously holds.
Q.e.d.
\section{Transitive and regular Courant algebroids}
In this section we determine the structure of the {\it transitive
Courant algebroids}, i.e., Courant algebroids with a surjective
anchor. The results may then be extended to {\it regular Courant
algebroids}, i.e., Courant algebroids with a constant rank anchor.
When this paper was ready, I was informed that the transitive
Courant algebroids were also determined by P. \v{S}evera in
unpublished correspondence with A. Weinstein, without metric
connections (see Remark \ref{obsextension} later on).

Let $(E,g,\rho)$ be a Courant vector bundle with a surjective
anchor $\rho$. Then, $K=ker\,\rho$ is a (regular) subbundle of $E$
of rank $k=r-n$, where $r=rank\, E$ and $n=dim\, M$, and, if
$K^{\perp_g}$ is the $g$-orthogonal subbundle of $K$, $dim(K\cap
K^{\perp_g})=r-s$, where $s=dim(K+K^{\perp_g})$. In view of the
properties of a Courant anchor $dim(im\,\partial)=n$ and
$im\,\partial\subseteq(K\cap K^{\perp_g})$ whence, it follows
easily that $s=k$. Accordingly, $K$ is a $g$-coisotropic subbundle
of $E$ and $K^{\perp_g}=im\,\partial\subseteq K$.

In this situation, it is known that there exists an isotropic,
complementary subbundle $Q$ of $K$ in $E$ and a complementary
subbundle $C$ of $im\,\partial$ in $K$ such that
\begin{equation}\label{Eincazcoiso}
E=P\oplus C\hspace{1cm}(P=im\,\partial\oplus Q)\end{equation} is a
$g$-orthogonal decomposition, the restriction of $g$ to $C$ is
non-degenerate, and the restriction of $g$ to $P$ is
non-degenerate and {\it neutral} (i.e., of signature zero). We
will say that $P$ is a {\it neutral completion of $im\,\partial$}.
(The reader may see
\cite{V87} for similar results in the symplectic case and the proofs
are the same in the pseudo-Euclidean case.) For what follows, we
fix a decomposition (\ref{Eincazcoiso}), which means that we also
have $E=K\oplus Q$, and $\rho|_Q$ is an isomorphism with the
inverse $\sigma: TM\rightarrow Q$. We will denote by
$p_K,p_Q,p_P,p_{im\,\partial},p_C$ the projections of $E$ onto the
corresponding subspaces, respectively.

Furthermore, we construct a metric connection $\nabla$ of $E$ as
follows. $\nabla$ will be the sum of metric connections of the
components $P,C$. Furthermore, the component $\nabla^P$ will be a
sum $\nabla^{im\,\partial}\oplus\nabla^Q$ where $\nabla^Q$ is
arbitrary and $\nabla^{im\,\partial}$ is defined by the condition
$$X(g(q,\partial f))=g(q,\nabla^{im\,\partial}_X(\partial
f)) + g(\nabla^Q_Xq,\partial f)\hspace{2mm}(q\in\Gamma
Q,X\in\Gamma TM,f\in C^\infty(M)).$$ This condition defines well
$\nabla^{im\,\partial}$ because $g|_P$ is neutral. A metric
connection on $(E,g)$ which is obtained by the process described
above is said to be {\it suitable}, and we fix one such suitable
connection.

The component $\nabla^Q$ of a suitable connection may be
identified with a linear connection $D$ on $M$ by means of the
formula
\begin{equation}\label{conexcorespD} D_XY=\rho(\nabla_X(\sigma Y))
\hspace{5mm}(X,Y\in\Gamma TM).
\end{equation}

The following formula defines a $3$-form
$B_1\in\Gamma\wedge^3E^*$:
\begin{equation}\label{Btranzitiv} B_1(e_1,e_2,e_3) =
\sum_{Cycl(1,2,3)}g(\sigma(T_{(\nabla,\rho)}(e_1,e_2)),p_{im\,\partial}(e_3)).
\end{equation} From the fact that the $\nabla$-parallel
translations preserve the subbundle $im\,\partial$, and formula
(\ref{pstorsiunepeE}), it follows that $B_1$ satisfies condition
(\ref{condicuB2}). Hence, formula (\ref{solgencubeta}) with
$\beta=\beta_1$ defined by $B_1$ yields a structure of a
pre-Courant algebroid on $(E,g,\rho)$, with a bracket that we
denote by $[\,,\,]_1$.

Like for any pre-Courant algebroid, in the transitive case too the
brackets $[e,\partial f]$ are always given by formula
(\ref{cptantisim}). Furthermore, we get
\begin{prop}\label{crosetQintranz} The bracket $[\,,\,]_1$
is defined by the formulas
\begin{equation}\label{prQptcoiso} p_Q[e_1,e_2]_1 = p_Q\nabla_{\rho
e_1}e_2-p_Q\nabla_{\rho e_2}e_1 -
\sigma(T_{(\nabla,\rho)}(e_1,e_2)),\hspace{5mm} \forall e_1,e_2\in
\Gamma E,\end{equation}
\begin{equation} \label{prKptcoiso} g(p_K[e_1,e_2]_1,e) =
g(\nabla_{\rho e_1}e_2,e) - g(\nabla_{\rho e_2}e_1,e)
\end{equation} $$- \frac{1}{2}g(\nabla_{\rho e}e_2,e_1) +
\frac{1}{2}g(\nabla_{\rho e}e_1,e_2)$$ $$- g(p_Q\nabla_{\rho
e_1}e_2,e) + g(p_Q\nabla_{\rho e_2}e_1,e)$$ $$+
g(\sigma(T_{(\nabla,\rho)}(e_1,e),e_2) -
g(\sigma(T_{(\nabla,\rho)}(e_2,e),e_1),$$ where $e\in\Gamma E$ and
we use an arbitrary decomposition {\rm(\ref{Eincazcoiso})} and an
arbitrary suitable connection $\nabla$.
\end{prop} \noindent {\bf Proof.}
The first formula follows by applying $\sigma$ to
(\ref{pstorsiunepeE}), since $\sigma\circ\rho$ is the projection
$p_Q$ and the bracket $[\,,\,]$ satisfies property i).

Furthermore, if we expend the expression of the bracket
$[\,,\,]_1$, take the scalar $g$-product by an arbitrary
$e\in\Gamma E$, and use formula (\ref{prQptcoiso}), we see that
(\ref{prKptcoiso}) also holds.

Formulas (\ref{prQptcoiso}) and (\ref{prKptcoiso}) define the
projections of the bracket $[\,,\,]_1$ on $Q$ and $K$, hence,
completely define the bracket. Q.e.d.

In the following proposition we give a more transparent expression
of the bracket $[\,,\,]_1$.
\begin{prop}\label{crosetfinalcaztranzitiv}
Let $(E,g,\rho)$ be a Courant vector bundle with a surjective
anchor, for which a choice of a spitting {\rm(\ref{Eincazcoiso})}
and of a suitable connection $\nabla$ is made. Then, $(E,g,\rho)$
has a structure of pre-Courant algebroid with the bracket
$[\,,\,]_1$ given by the formulas
\begin{equation}\label{finalbracket}
\begin{array}{l}[q_1+\partial f_1,q_2+\partial
f_2]_1=\sigma[\rho q_1,\rho q_2]+\frac{1}{2}[\partial((\rho
q_1)f_2) -\partial((\rho q_2)f_1)],
\vspace{2mm}\\

[c_1,c_2]_1=-\gamma(c_1,c_2),\hspace{5mm} [c,q+\partial
f]_1=-\nabla_{\rho q}c, \end{array}
\end{equation} where $f,f_1,f_2\in C^\infty(M);\,q,q_1,q_2\in\Gamma
Q;\, c,c_1,c_2\in\Gamma C$ and $\gamma$ is defined by
{\rm(\ref{defluigamma})}. Furthermore, $[\,,\,]_1$ is a Courant
algebroid bracket iff the connection $\nabla$ is flat.
\end{prop}
\noindent {\bf Proof.} In view of the
properties of decomposition (\ref{Eincazcoiso}), $k\in K$ is
completely defined by the scalar products $g(k,q)$ $\forall q\in
Q$. But, if $e,e_1,e_2$ of (\ref{prKptcoiso}) are in $Q$, using
the definition of a suitable connection we get $p_K[q_1,q_2]_1=0$.

Accordingly, with (\ref{prQptcoiso}) and (\ref{pstorsiunepeE}) we
get $[q_1,q_2]_1=\sigma[\rho q_1,\rho q_2]$ and the first formula
(\ref{finalbracket}) follows if we also take into account
(\ref{cptantisim}).

By similar considerations based on the properties of decomposition
(\ref{Eincazcoiso}) and of a suitable connection,
(\ref{prQptcoiso}), (\ref{prKptcoiso}) and (\ref{cptantisim})
yield the remaining formulas (\ref{finalbracket}).

Of course, brackets $[\,,\,]_1$ with general factors
$\sum_ih_i\partial f_i\in im\,\partial$ will be deduced from
(\ref{finalbracket}) by means of property iv), Proposition
\ref{propLiu}. (Alternatively, we may use again (\ref{prQptcoiso})
and (\ref{prKptcoiso}).) Notice also that (\ref{defluigamma})
implies $\gamma(e_1,e_2)\in K^{\perp_g}=im\,\partial$, $\forall
e_1,e_2\in\Gamma E$.

From the first formula (\ref{finalbracket}), we see that
$(P,g|_{P},
\rho|_{P})$ $(P=Q\oplus im\,\partial)$ is a
pre-Courant algebroid with the induced bracket. Moreover, it is
easy to check that property iii) of Proposition \ref{propLiu} also
holds on $P$, therefore, we actually have a Courant algebroid $P$,
and (\ref{finalbracket}) is the structure defined on $P\oplus Q$
by Proposition \ref{sugerat2}. Accordingly, the last assertion of
the present Proposition follows from the last assertion of
Proposition \ref{sugerat2}. Q.e.d.

As shown by Proposition \ref{strsugerate}, all the other
pre-Courant brackets of the Courant bundle with surjective anchor
$(E,g,\rho)$ will be obtained from (\ref{finalbracket}) by the
addition of a form $\lambda$ that satisfies the corresponding
hypotheses, whence, $\lambda\in\Gamma(\wedge^2E^*\otimes K)$ and
(\ref{condilambda}) holds. The corresponding formulas are
\begin{equation}\label{finalbracket1}
\begin{array}{l}[q_1+\partial f_1,q_2+\partial
f_2]=\sigma[\rho q_1,\rho q_2]+\lambda(q_1,q_2)\vspace{2mm}\\

\hspace{35mm}+\frac{1}{2}[\partial((\rho q_1)f_2)
-\partial((\rho q_2)f_1)],
\vspace{2mm}\\

[c_1,c_2]=-\gamma(c_1,c_2)+\lambda(c_1,c_2),\hspace{5mm}
[c,q+\partial f]=-\nabla_{\rho q}c+\lambda(c,q). \end{array}
\end{equation}

In order to get Courant algebroid brackets we must ask $\lambda$
to satisfy condition (\ref{condcociclu0}). If expressed on
arguments $q\in\Gamma Q, c\in\Gamma C$, (\ref{condcociclu0})
decomposes into the following four components
\begin{equation}\label{component1}
(\partial_{[\,]}\Lambda)(q_1,q_2,q_3)=0,\end{equation}
\begin{equation}\label{component2}
(\partial_{[\,]}\Lambda)(q_1,q_2,c)=-R_\nabla(\rho q_1,\rho q_2)c,
\end{equation}
\begin{equation}\label{component3}
(\partial_{[\,]}\Lambda)(q,c_1,c_2)=
\partial[g(q,\gamma(c_1,c_2))]_1 + [q,\gamma(c_1,c_2)]_1\end{equation}
$$-\gamma(\nabla_{\rho q}c_1,c_2) - \gamma(c_1,\nabla_{\rho
q}c_2),$$
\begin{equation}\label{component4}
(\partial_{[\,]}\Lambda)(c_1,c_2,c_3)=
\sum_{Cycl(1,2,3)}[c_1,\gamma(c_2,c_3)]_1.\end{equation}

In principle, formulas (\ref{finalbracket1}) and
(\ref{component1})-(\ref{component4}) yield all the transitive
Courant algebroids.

We get a more transparent result if we restrict ourselves to the
subclass of transitive Courant algebroids that admit a
bracket-closed neutral extension of the subbundle $im\,\partial$.
We shall call them {\it transitive, restricted, Courant
algebroids}.

The general bracket (\ref{finalbracket1}) is restricted iff it is
defined by a form $\lambda$ such that, $\forall q_1,q_2\in\Gamma
Q$, $\lambda(q_1,q_2)\in im\,\partial$. This condition is
equivalent with $\Lambda(q_1,q_2,c)=0$, which, because of the
skew-symmetry, is equivalent with $\lambda(c,q)\in
\Gamma C$. Then, we may change the $C$-component of the connection
$\nabla$ by $\nabla_XC\mapsto \nabla_XC+\lambda(c,\sigma(X))$ and
get for the same bracket a simplified expression
(\ref{finalbracket1}) that looks as if we have used an additional
form $\lambda$ such that $\lambda(c,q)=0$ and, accordingly,
$\lambda(c_1,c_2)\in\Gamma C$. Then, conditions
(\ref{component1})-(\ref{component4}) become simpler and we obtain
\begin{prop}\label{finalptalgCtranz} The bracket
{\rm(\ref{finalbracket1})} defines a transitive, restricted,
Courant algebroid iff the $C$-component of the connection $\nabla$
is flat and, in addition to the conditions for a pre-Courant
algebroid, the form $\lambda$ also satisfies the conditions
\begin{equation}\label{condrestrict}
\lambda(c,q)=0,\;\lambda(q_1,q_2)\in\Gamma(im\,\partial),\;
\lambda(c_1,c_2)\in\Gamma C,\end{equation}
\begin{equation}\label{condcociclu}
\partial(\Lambda(q_1,q_2,q_3))=
\sum_{Cycl(1,2,3)}\{ \lambda(\sigma[\rho q_1,\rho q_2],q_3) -
[q_3,\lambda(q_1,q_2)]\},
\end{equation}
\begin{equation} \label{condcociclu2} \partial(\Lambda(c_1,c_2,c_3))=
-\sum_{Cycl(1,2,3)}\gamma(\lambda(c_1,c_2),c_3), \end{equation}
\begin{equation} \label{condcociclu3} (\nabla_{\rho
q}\lambda)(c_1,c_2)=0.\end{equation}
\end{prop} \noindent {\bf Proof.}
Conditions (\ref{condrestrict}) transform (\ref{component2}) into
the flatness of $\nabla$, (\ref{component1}) into
(\ref{condcociclu}), (\ref{component4}) into (\ref{condcociclu2})
and (\ref{component3}) into (\ref{condcociclu3}). Q.e.d.
\begin{corol}\label{corolarC} A Courant vector bundle with
surjective anchor has a restricted Courant algebroid structure iff
it is a Whitney sum {\rm(\ref{Eincazcoiso})} of a Courant
subalgebroid with neutral metric and a flat pseudo-Euclidean
bundle.
\end{corol} \noindent {\bf Proof.} For a Courant algebroid, the
conditions stated by the corollary were proven in Proposition
\ref{finalptalgCtranz}. The bracket of the subalgebroid
is defined by the first formula (\ref{finalbracket1}). Conversely,
if the conditions hold, a corresponding Courant algebroid bracket
is defined by (\ref{finalbracket1}) with $\lambda=0$. Q.e.d.
\begin{rem}\label{rangminimal} {\rm The minimal possible rank of a
transitive Courant algebroid over an $n$-dimensional manifold $M$
is $2n$. Proposition \ref{finalptalgCtranz} yields all the
transitive Courant algebroids of minimal rank. Namely, they have
no $C$-component and the bracket is defined by the first formula
(\ref{finalbracket1}) where $\lambda(q_1,q_2)\in im\,\partial$
satisfies the condition (\ref{condcociclu}).}\end{rem}
\begin{rem}\label{obsextension} {\rm With the notation used above,
let $E$ be a transitive Courant algebroid. Then $E/im\,\partial$
gets an induced structure of a transitive Lie algebroid and we
have the exact sequence
$$0\rightarrow im\,\partial\stackrel{\subseteq}{\rightarrow}E
\stackrel{\pi}{\rightarrow}E/im\,\partial\rightarrow0.$$ Hence $E$
may be seen as a central extension of a transitive Lie algebroid.
In the restricted case, the splitting $\sigma$ induces a splitting
$TM\stackrel{\sigma}{\rightarrow}E/im\,\partial$ which is a
morphism of Lie algebroids. P. \v{S}evera expressed the bracket of
the Courant algebroid $E$ by means of the bracket of the Lie
algebroid $E/im\,\partial$.}\end{rem}

The results obtained so far in this section straightforwardly
extend to regular Courant algebroids. Indeed, such an algebroid
$(E,g,\rho)$ is a transitive Courant algebroid over the base
manifold $M$ of $E$ seen as the sum of the leaves of the foliation
$\mathcal{F}=im\,\rho$. Therefore, we get formulas
(\ref{finalbracket1}) again. In order to ensure that we obtain
brackets that are differentiable with respect to the original
differentiable structure of $M$ it suffices to use metric
connections and forms $\lambda,\Lambda$ that enjoy this kind of
differentiability. Such connections are just Lie algebroid
connections for the tangent Lie algebroid of $\mathcal{F}$ (called
$\mathcal{F}$-partial connections or connections along the leaves
of $\mathcal{F}$ in foliation theory). Thus, Propositions
\ref{crosetfinalcaztranzitiv},
\ref{finalptalgCtranz} and Corollary \ref{corolarC}, where we ask
the anchor to be surjective over a regular foliation $\mathcal{F}$
of $M$ and the connection to be along the leaves of $\mathcal{F}$,
describe all the regular (restricted) Courant algebroids.\\
\section{Miscellanies}
{\it a)} The basic example of a skew-symmetric Courant algebroid
appeared in \cite{C}. It was the vector bundle $E=TM\oplus T^*M$,
endowed with the neutral pseudo-Euclidean metric
\begin{equation}\label{gCourant}
g(X\oplus\alpha,Y\oplus\beta) =\frac{1}{2}(i(X)\beta+i(Y)\alpha),
\end{equation} the non-degenerate cross section
$\omega\in\wedge^2E^*$,
\begin{equation}\label{omegaCourant}
\omega(X\oplus\alpha,Y\oplus\beta)
=\frac{1}{2}(i(X)\beta-i(Y)\alpha), \end{equation} the {\it
Courant bracket}
\begin{equation}\label{crosetC}
[X\oplus\alpha,Y\oplus\beta]=[X,Y]\oplus[L_X\beta-L_Y\alpha
-d(\omega(X\oplus\alpha,Y\oplus\beta))],\end{equation} and the
projection $\rho(X\oplus\alpha)=X$, where $X,Y\in
\Gamma(TM),\,\alpha,\beta\in\Gamma(T^*M)$ and, for clarity, we
denoted an element $(X,\alpha)=X+\alpha\in TM\oplus T^*M$ by
$X\oplus\alpha$.

Straightforward computations show that the bracket
(\ref{crosetC}), together with $g,\rho$, satisfy the conditions of
Definition \ref{skewCalg}, and $TM\oplus T^*M$ is a skew-symmetric
Courant algebroid.

For the same data, the formula
\begin{equation}\label{brCnonskew} (X\oplus\alpha)\star
(Y\oplus\beta) = [X,Y]\oplus(L_X\beta-i(Y)d\alpha)\end{equation}
defines a structure of a non-skew-symmetric Courant algebroid and
the Courant bracket (\ref{crosetC}) is the skew-symmetric part of
the bracket (\ref{brCnonskew}).

The Courant bracket was extended by \v Severa and Weinstein
\cite{SW} by the addition of a term of the form $i(X\wedge Y)\Phi$
to the $T^*M$-component of the right hand side of (\ref{crosetC}),
$\Phi$ being a closed $3$-form on $M$.

The Courant algebroid structures described above are transitive
hence, particular cases of the general formulas
(\ref{finalbracket}), (\ref{finalbracket1}). Indeed, in the
present case, $\partial f=0\oplus df$ $(f\in C^\infty(M))$,
$im\,\partial=T^*M$, $C=\{0\}$, we may take $Q=TM$, and a suitable
connection is provided by any linear connection on $M$. Then, the
first formula (\ref{finalbracket}) becomes (\ref{crosetC}).
Indeed, this is trivial for $\alpha=df_1,\beta=df_2$ $(f_1,f_2\in
C^\infty(M))$, and it is true for any $\alpha,\beta$ because the
two formulas (\ref{finalbracket}) and (\ref{crosetC}) behave in
the same way when arguments are multiplied by a function. The
addition of a form $\lambda$ leads to the \v Severa-Weinstein
Courant bracket with $\Phi=-2\Lambda|_Q$. Since $C=\{0\}$, we
necessarily are in the restricted case, $\lambda,\Lambda$ satisfy
(\ref{Lambda}), and the conditions
(\ref{condcociclu})-(\ref{condcociclu3}) reduce to $d\Phi=0$.
Therefore, the Courant and \v Severa-Weinstein brackets define all
the Courant algebroid structures on $TM\oplus T^*M$ endowed with
the metric (\ref{gCourant}) and the anchor
$\rho(X\oplus\alpha)=X$.

However, $TM\oplus T^*M$ may have more Courant algebroid
structures if, for instance, we change the anchor, as we will see
below.

In \cite{Liu}, the Courant bracket was extended to vector bundles
$A\oplus A^*$, where $(A,A^*)$ is a Lie bialgebroid  with anchors
$\alpha,\alpha^*$, respectively, such that the extended bracket,
the metric $g$ defined like in (\ref{gCourant}) and the anchor
$\rho=\alpha+\alpha^*$ define a structure of skew-symmetric
Courant algebroid. The extended bracket is
\begin{equation} \label{eqXu1}
[\hspace{-1pt}[a\oplus a^*,b\oplus
b^*]\hspace{-1pt}]=\{[a,b]_A+L^{*}_{a^*}b -
L^{*}_{b^*}a+d^{*}(\omega(a\oplus a^*,b\oplus
b^*))\}\end{equation}
$$\oplus \{[a^*,b^*]_{A^*}+L_{a}b^* -
L_ba^{*}-d(\omega(a\oplus a^*,b\oplus b^*))\},$$ where $d,L$,
respectively $d^*,L^*$, are the exterior differential and Lie
derivative associated with the Lie algebroid structure of $A$,
respectively $A^*$, and $\omega$ is defined like in
(\ref{omegaCourant}).

One can check that the bracket (\ref{eqXu1}) is the skew-symmetric
part of the product
\begin{equation} \label{eqXu10}
(a\oplus a^*)\star(b\oplus b^*)=\{[a,b]_A+L^{*}_{a^*}b +
i(b^*)d^*a\}\end{equation}
$$\oplus \{[a^*,b^*]_{A^*}+L_{a}b^* - i(b)da^*\}.$$

It is possible to connect the bracket (\ref{eqXu1}) with a metric
connection as indicated in Proposition \ref{crosetcuA}, but this
doesn't seem to give interesting formulas. On the other hand, we
may use the bracket (\ref{eqXu1}) and Propositions
\ref{strsugerate}, \ref{strsugerate0} in order to derive new
Courant algebroid structures on the vector bundles $A\oplus A^*$
and $A\oplus A^*\oplus C$ for any flat pseudo-Euclidean bundle
$C$.

The Courant bracket (\ref{crosetC}) on $TM\oplus TM^*$ is the
particular case of the bracket (\ref{eqXu1}) where $A=TM$ with the
Lie bracket and $A^*=T^*M$ with the zero bracket and zero anchor.
But, if we assume that $P$ is a Poisson bivector field on $M$ and
use the cotangent algebroid structure defined by $P$ on $T^*M$, we
get a new Courant algebroid structure on $TM\oplus TM^*$ with the
same metric (\ref{gCourant}), with the anchor
$Id_{TM}\oplus\sharp_P$, and with the bracket (\ref{eqXu1}).
Again, we are in the transitive, restricted case, and it will be
possible to express the bracket under the form
(\ref{finalbracket}). Indeed, the kernel of the anchor is
$$K=\{X\oplus\alpha\,/\,X+\sharp_P\alpha=0\}\approx T^*M,$$
and we may use the complementary subbundle $Q=\{X\oplus0\}\approx
TM$ such that a cross section of $TM\oplus TM^*$ decomposes as
follows
$$X\oplus\alpha=((X+\sharp_P\alpha)\oplus0)+
((-\sharp_P\alpha)\oplus\alpha).$$ It is easy to see that
$$\partial f=(-X^P_f)\oplus df,\hspace{1cm}f\in C^\infty(M),$$
where $X^P_f$ is the Hamiltonian vector field of $f$ with respect
to $P$. Then, by technical calculations, one checks that the
brackets defined by the first formula (\ref{finalbracket}) and by
(\ref{eqXu1}) coincide; the checks are to be made in each of the
cases: two arguments in $Q$, one in $Q$ and one in $im\,\partial$
and two arguments in $im\,\partial$.

The original Courant bracket also leads to an example of a
non-restrictive, transitive pre-Courant algebroid. Assume that the
manifold $M$ is endowed with a Riemannian metric $G$ and consider
the vector bundle $E=TM\oplus T^*M\oplus TM$ with the anchor
defined as the projection on the first term and the metric
$g\oplus G$, where $g$ is given by (\ref{gCourant}). Then,
$\partial f=0\oplus df\oplus0$ and $E$ is a Courant vector bundle
with surjective anchor, with the natural decomposition $E=(Q\oplus
im\,\partial)\oplus C$ where $Q=TM$ and $C=(TM,G)$, and the
suitable connection defined by the Levi-Civita connection $\nabla$
of $G$, which satisfies $T_{\nabla,\rho}=0$.

Accordingly, we get a pre-Courant algebroid bracket on $TM\oplus
T^*M\oplus TM$ if we use the corresponding formulas
(\ref{finalbracket}). The first formula (\ref{finalbracket}) again
yields the original Courant bracket and the remaining brackets
(\ref{finalbracket}) are determined by the Levi-Civita connection
$\nabla$ of $G$ and the value of $\gamma$ as defined by
(\ref{defluigamma}), which yields
\begin{equation}\label{gammaptG} \gamma(X_1\oplus\alpha_1\oplus
Y_1,X_2\oplus\alpha_2\oplus Y_2)=0\oplus\xi(Y_1,Y_2)\oplus0,
\end{equation} where $$\xi(Y_1,Y_2)=G(Y_1,\nabla Y_2) - G(\nabla
Y_1,Y_2)$$ and $X_a,Y_a\in\Gamma TM,\,\alpha_a\in\Gamma T^*M$
$(a=1,2)$.

Now, let $\Phi$ be a differential $3$-form on $M$ and define
\begin{equation}\label{lambdaptG} \lambda(X_1\oplus\alpha_1\oplus
Y_1,X_2\oplus\alpha_2\oplus Y_2) \end{equation}
$$=0\oplus2[i(X_1\wedge Y_2-X_2\wedge Y_1)\Phi]
\oplus\sharp_G[i(X_1\wedge X_2)\Phi].$$ This form vanishes if one
of the arguments belongs to $im\,\partial$ and
\begin{equation}\label{LambdaptG}
\Lambda(X_1\oplus\alpha_1\oplus Y_1,X_2\oplus\alpha_2\oplus
Y_2,X_3\oplus\alpha_3\oplus Y_3)\end{equation}
$$=\Phi(X_1,X_2,Y_3)+\Phi(X_1,Y_2,X_3)+\Phi(Y_1,X_2,X_3)$$
is skew symmetric. Therefore $\lambda$ of (\ref{lambdaptG}) may
serve as an additional term that leads to a new pre-Courant
algebroid bracket (\ref{finalbracket1}) on $TM\oplus T^*M\oplus
TM$, which is not restricted. It will be a Courant algebroid
bracket iff the conditions (\ref{component1})-(\ref{component4})
hold.
The significance of these conditions is unclear.\\

{\it b)} The Courant bracket (\ref{crosetC}) may be defined for
every Lie algebroid $A$ and it is the particular case of
(\ref{eqXu1}) where the dual bundle $A^*$ is endowed with the zero
anchor and the zero bracket. In what follows, we use this remark
in order to define a Courant algebroid structure on the tangent
bundle of a {\it para-Hermitian manifold} $M$. We recall the
definition \cite{CFG}: a {\it para-Hermitian structure} on $M$
consists of a neutral metric $g$ on $TM$ and a $(1,1)$-tensor
field $F$ that satisfies the conditions
\begin{equation}\label{eqparaherm} F^2=I=Id.,\;
g(FX,FY)=-g(X,Y),\hspace{5mm}
\forall X,Y\in \Gamma TM,\end{equation}
\begin{equation}\label{Nij}
N_F(X,Y)=[FX,FY]-F[FX,Y]-F[X,FY]+F^2[X,Y]=0.\end{equation}
Condition (\ref{Nij}), which is the vanishing of the {\it
Nijenhuis tensor} of $F$, is the integrability condition of the
structure.

It follows that $M$ also has a non degenerate $2$-form
\begin{equation} \label{2-forma} \omega(X,Y)=g(FX,Y),\end{equation}
which satisfies the condition
\begin{equation}\label{eqparaherm2} \omega(FX,FY)=-\omega(X,Y),\hspace{5mm}
\end{equation} and that $TM=W_{+}\oplus W_{-}$,
where the terms are the $\pm1$-eigendistributions of $F$ and are
integrable because of (\ref{Nij}). We will denote by
\begin{equation} \label{projectors} F_{\pm}=\frac{1}{2}(I\pm F)
\end{equation} the projectors on $W_{\pm}$, respectively.
From (\ref{eqparaherm}) and (\ref{eqparaherm2}) it follows that
$W_{\pm}$ are maximal isotropic subbundles with respect to $g$ and
Lagrangian subbundles with respect to $\omega$. Accordingly, the
musical isomorphism $\flat_g$ sends $W_{\pm}$ onto the dual space
$W_{\mp}^*$ and defines an isomorphism $TM\approx W_{+}\oplus
W^*_{+}$.

Because of integrability, $W_{+}$, with the Lie bracket, is a Lie
algebroid and, if we use (\ref{crosetC}) in this case, we get a
bracket on the tangent bundle $TM$ defined by
\begin{equation}\label{crosetCcuF} [X,Y]_{F_+}=[X_{+},Y_{+}]
+\sharp_g\{L_{X_{+}}\flat_gY_{-}-L_{Y_{+}}\flat_gX_{-}-
\frac{1}{2}d[\omega(X,Y)]\},
\end{equation} where $X_{\pm}=F_{\pm}X,\,Y_{\pm}=F_{\pm}Y$.
The conclusion is that $(TM,g,F_+,[\,,\,]_{F_+})$ is a regular,
skew-symmetric Courant algebroid.
\begin{rem}\label{FptCourant} {\rm In an obvious way, we may speak
of {\it para-Hermitian vector bundles} and, for any differentiable
manifold $M$, the bundle $TM\oplus T^*M$ is para-Hermitian with
the metric (\ref{gCourant}) and with $F$ defined by
\begin{equation}\label{FCourant}
F(X\oplus\alpha)=X\oplus(-\alpha).
\end{equation} Then, formulas (\ref{crosetC}) and
(\ref{crosetCcuF}) are similar.}\end{rem}

{\it c)} In what follows we give some indications about the use of
Proposition
\ref{crosetcuA} for possible Courant algebroid structures
on the tangent bundle of an arbitrary differentiable manifold $M$.

A Courant algebroid structure on $TM$ requires a pseudo-Riemannian
metric $g$ on $M$, and a field $\phi$ of endomorphisms of $TM$,
which will be the anchor. Notice that $\partial
f=(1/2)\sharp_g(df\circ\phi)$ and that
$span\{df\circ\phi\}=ann\,ker\,\phi$. Hence, $\phi$ is a Courant
anchor iff $im\,\phi$ is a generalized, completely integrable
distribution and $\sharp_g(ann\,ker\,\phi)$ is an isotropic
distribution. In particular, we must have $rank\,\phi\leq
(1/2)dim\,M$. Furthermore, one must have a bracket
$[\hspace{-1pt}[X,Y]\hspace{-1pt}]$
$(X,Y,[\hspace{-1pt}[X,Y]\hspace{-1pt}] \in\Gamma TM)$ that
satisfies properties i),iii), v) of Proposition \ref{propLiu}. In
Section 4{\it b)} we had an example of such a situation.

Generally, on $TM$ we have the Levi-Civita connection $\nabla$ of
$g$ , which is metric, and we may use it to express the bracket.
The connection $\nabla$ has no torsion but it has a $\phi$-torsion
\begin{equation}\label{pseudoTpeTM} \begin{array}{lcl}
T_{(\nabla,\phi)}(X,Y)&=&\phi(\nabla_{\phi X}Y-\nabla_{\phi
Y}X)-[\phi X,\phi Y]\vspace{2mm}\\&=&(\nabla_{\phi
Y}\phi)(X)-(\nabla_{\phi X}\phi)(Y).\end{array}\end{equation} A
technical calculation shows that the $\phi$-torsion is related
with the Nijenhuis tensor $N_\phi$ (see (\ref{Nij}). Namely,
\begin{equation}\label{TcuN} T_{(\nabla,\phi)}(X,Y) =
(\phi\circ\nabla_Y\phi)(X) - (\phi\circ\nabla_X\phi)(Y) -
N_\phi(X,Y).\end{equation}

The operator $\gamma$ of (\ref{defluigamma}) will now be defined
by
\begin{equation}\label{gammacuLC} g(\gamma(X,Y),Z) =
\frac{1}{2}[g(X,\nabla_{\phi Z}Y)-g(\nabla_{\phi
Z}X,Y)],\end{equation} and the Courant algebroid bracket will be
\begin{equation}\label{CcuLC}
[\hspace{-1pt}[X,Y]\hspace{-1pt}]= \nabla_{\phi X}Y- \nabla_{\phi
Y}X - \gamma(X,Y) -\beta(X,Y).
\end{equation} where $\beta$ is a tensor field of type (1,2) on
$M$ such that $\phi(\beta(X,Y))= T_{(\nabla,\phi)}(X,Y)$ and
$B(X,Y,Z)=g(\beta(X,Y),Z)$ is a $3$-form on $M$.

Thus, essentially, a Courant algebroid structure on $(TM,g,\phi)$
is a special type of a $3$-form. What we must still ask is the
fulfillment of condition $\mathcal{C}=0$ for $\mathcal{C}$ defined
by (\ref{CprinCzero}). Unfortunately, this condition is too
complicated and does not provide a practical way to find new
Courant algebroids. This is true even if stronger conditions are
added. For instance, if we ask $\phi$ to be $\nabla$-parallel
along paths in the leaves of $im\,\phi$, (\ref{pseudoTpeTM}) shows
that $T_{(\nabla,\phi)}=0$ and we may try the solution $\beta=0$.
The resulting bracket (\ref{CcuLC}) defines a Courant algebroid
structure iff (see (\ref{Czero}))
$$\sum_{Cycl(X,Y,Z)}\{R_\nabla(\phi X,\phi Y)Z - (\nabla_{\phi
Z}\gamma)(X,Y) - \gamma(\gamma(X,Y),Z)\}=0.$$

{\it d)} Now, in a different direction, we refer to the case of
what should be called {\it foliated, regular Courant algebroids},
and show another way of describing their brackets.

Let $M$ be a manifold endowed with a regular foliation
$\mathcal{F}$. Let $E$ be a foliated, vector bundle over $M$
endowed with a foliated, pseudo-Euclidean metric $g$. This means
that $E$ has a given, maximal system of local trivializations with
transition functions that are constant along the leaves of
$\mathcal{F}$ (foliated functions), and, $\forall
e_1,e_2\in\Gamma_{pr}E$, where $\Gamma_{pr}$ denotes the space of
projectable cross sections of $E$ (i.e., constant along the
leaves), the function $g(e_1,e_2)$ is constant along the leaves.
(E.g., see \cite{Mol} for the theory of foliations.) Then, we have
\begin{prop}\label{exitconexptfol} On  a foliated pair $(E,g)$
there exist metric connections $\nabla$  that satisfy the
condition $\nabla_Xe=0$, $\forall X\in\Gamma T\mathcal{F}$ and
$\forall e\in\Gamma_{pr}E$.\end{prop} \noindent {\bf Proof.} Take
an arbitrary normal bundle $N\mathcal{F}$ (i.e.,
$TM=T\mathcal{F}\oplus N\mathcal{F}$) and an arbitrary metric
connection $\tilde\nabla$ of $(E,g)$. Then define
\begin{equation}\label{defconexptfol}
\nabla_Xe=\left\{\begin{array}{ll}\tilde\nabla_Xe& \forall X\in
\Gamma N\mathcal{F},\;\forall e\in\Gamma E \vspace{2mm}\\ 0&
\forall X\in \Gamma T\mathcal{F},\;\forall
e\in\Gamma_{pr}E.\end{array} \right. \end{equation} It is easy to
check that this produces a connection as required. Indeed, we must
still define $\nabla_Xe$ for $X\in\Gamma T{\mathcal F}$ and a non
projectable $e\in\Gamma E$. For this purpose, we take local
projectable bases $e_\alpha$ of $E$ and, for
$$e=\sum_\alpha f_\alpha e_\alpha,\hspace{1cm}
e_\alpha\in\Gamma_{pr}E,\;f_\alpha\in C^\infty(M),$$ we put
$$\nabla_Xe=
\sum_\alpha(Xf_\alpha)e_\alpha.$$ Q.e.d.

A connection that satisfies the properties stated by Proposition
\ref{exitconexptfol} will be called an {\it adapted connection} of
$(E,g)$.

Let $(E,g)$ be foliated and let us assume that there exists a
surjective morphism $\rho:E
\rightarrow T\mathcal{F}$, which is a Courant anchor of $(E,g)$.
Then, for any adapted connection $\nabla$, one has a non-zero
$\rho$-torsion given by
\begin{equation}\label{torsadaptata} T_{(\nabla,\rho)}(e_1,e_2)=
-[\rho e_1,\rho e_2],\hspace{1cm} e_1,e_2\in\Gamma_{pr}E.
\end{equation}
Formula (\ref{defluigamma}) yields the operator
\begin{equation}\label{cazulgammazero}
\gamma(e_1,e_2)=0,\hspace{1cm}
e_1,e_2\in\Gamma_{pr}E,\end{equation} and (\ref{solparticcugamma})
gives a bracket such that
\begin{equation}\label{cazparticularzero}[e_1,e_2]_0=0,\hspace{1cm}
e_1,e_2\in\Gamma_{pr}E.\end{equation} The values of this bracket
for arbitrary cross sections $f_1e_1,f_2e_2$,
$e_1,e_2\in\Gamma_{pr}E, f_1,f_2\in C^\infty(M)$ follows from
property iv), Proposition
\ref{propLiu}. Furthermore, the operator $\gamma$ satisfies the
condition $\mathcal{C}_0=0$, where $\mathcal{C}_0$ is defined by
(\ref{Czero}); this is obvious for projectable arguments and is
true for arbitrary arguments because $\mathcal{C}_0$ is a tensor.
Hence, we get
\begin{prop}\label{CourantpeEfol} For a triple $(E,g,\rho)$ as
described above, a Courant algebroid bracket is a bracket of the
form (\ref{solgencubeta}), where
$\beta\in\Gamma(\wedge^2E^*\otimes E)$ is associated with a
$3$-form $B\in\Gamma(\wedge^3E^*)$ and satisfies the conditions
\begin{equation}\label{condptBinEfol1}
\rho(\beta(e_1,e_2))=-[\rho e_1,\rho e_2]\end{equation}
\begin{equation}\label{condptBinEfol2}
\sum_{Cycl(1,2,3)}\beta(\beta(e_1,e_2),e_3)
=\sum_{Cycl(1,2,3)}[\beta(e_1,e_2),e_3]_0,
\end{equation} $\forall e_1,e_2,e_3\in\Gamma_{pr}E$.
\end{prop}
\section{Appendix: Dirac linear spaces}
The Courant algebroids resulted from the process of studying {\it
Dirac structures}, which are a significant generalization of the
Poisson structures
\cite{C}. Although this is not a subject of the present paper, we
have added this appendix, which shows that the known linear
algebra of Dirac structures is a part of para-Hermitian linear
algebra.

A para-Hermitian vector space is a $2n$-dimensional vector space
$W$ that has the structure indicated for the fibers of the tangent
bundle of a para-Hermitian manifold in the previous section. On
$W$ we have the ingredients $g,\omega,F,W_{\pm},F_{\pm}$, with the
algebraic properties stated in Section 4{\it b)}, and
$W=W_{+}\oplus W_{-}$, $\flat_g:W_{\pm}\approx W^*_{\mp}$.

The space $W$ has {\it adapted bases} $(b_i,c_j)$ $(i,j=1,...,n)$,
where $(b_1,...,b_n)$ is a basis of $W_{+}$, therefore, $\flat_g
b_1,...,\flat_g b_n$ is a basis of $W_{-}^*$, and $(c_1,...,c_n)$
is the corresponding dual basis of $W_{-}$, i.e.,
$g(b_i,c_j)=\delta_{ij}$ $(i,j=1,...,n)$.
\begin{defin} \label{defDirac} {\rm 1) \cite{C} A maximal, $g$-isotropic
subspace $L$ of the para-Hermiti\-an vector space $W$ is called a
{\it Dirac subspace} of $W$. 2) \cite{JR} A pair $(L,L')$ of
complementary Dirac subspaces of $W$ $(W=L\oplus L')$ is a {\it
reflector} in $W$.}
\end{defin}
\begin{prop} \label{Dtransvers} For any Dirac subspace
$L\subseteq W$, there exists a family $\mathcal{T}(L)$ of Dirac
subspaces that are complementary to $L$ in $W$, and
$\mathcal{T}(L)$ is an affine space modelled over the linear space
of the skew-symmetric matrices of order $n$. \end{prop}
\noindent {\bf Proof.} The results are analogous to known results
for Lagrangian subspaces of a symplectic vector space, and we
prove them as in the latter case, e.g., \cite{V87}. For any
subspace $S\subseteq W$ such that $W=L\oplus S$ and any basis
$(l_1,...,l_n)$ of $L$, there exists a unique {\it conjugated
basis} $(s_1,...,s_n)$ of $S$, such that $g(l_i,s_j)=\delta_{ij}$.
Using the conjugated basis, we can obtain vectors
$$u_i=s_i+\tau_i^kl_k$$ (the Einstein summation convention holds)
such that $g(u_i,u_j)=0$, and these vectors span a Dirac subspace
$L'$ that satisfies $L\oplus L'=W$. Furthermore, if $L''$ is
another Dirac subspace such that $L\oplus L''=W$ and if $(u_i)$ is
the conjugated basis of $(l_i)$ in $L'$ and $(v_i)$ is the
conjugated basis of $(l_i)$ in $L''$, there exists a unique
skew-symmetric matrix $(\theta_i^j)$ such that
$v_i=u_i+\theta_i^jl_j$. Q.e.d.

Using adapted bases of $W$, it follows that the set of the
reflectors of $W$ is the $n(n-1)$-dimensional homogeneous space
$\mathcal{R}=O(W,g)/pH(W)\approx O(n,n)/Gl(n)$, where $O(W,g)$ is
the $g$-preserving subgroup of the general linear group $Gl(W)$,
which acts transitively on $\mathcal{R}$, $pH(W)$ is the {\it
para-Hermitian subgroup}, which commutes with $F$ and is the
isotropy subgroup of the pair $(W_{+},W{-})\in\mathcal{R}$,
$O(n,n)\approx O(W,g)$ is the subgroup of $Gl(2n)$ which preserves
the canonical neutral metric, and $Gl(n)\approx pH(W)$ by the
embedding
\begin{equation}
\label{groupn} A\mapsto \left(
\begin{array}{cc}A&0\vspace{1mm}\\0&^tA^{-1}
\end{array} \right)\hspace{5mm}(A\in
Gl(n)) \end{equation}  \cite{{JR},{V98}}.

Proposition \ref{Dtransvers} shows that the set $\mathcal{D}$ of
Dirac subspaces of $W$ is the quotient space of $\mathcal{R}$ by
the equivalence relation with equivalence classes
$\mathcal{T}(L)$, hence, $\mathcal{D}$ is a
$[n(n-1)/]2$-dimensional space namely, the homogeneous space
$O(W,g)/O_{W_{+}}(W,g)$, where the isotropy group of
$W_{+}\in\mathcal{D}$ at the denominator is that of the elements
$\phi\in O(W,g)$ which satisfy the condition $F_{-}\circ\phi\circ
F_{+}=0$.

We also notice that $L\subseteq W$ is a Dirac subspace iff $F(L)$
is the $\omega$-orthogonal subspace of $L$. Therefore, if $L$ is a
Dirac subspace, $ker(\omega|_L)=L\cap F(L)$. On the other hand, it
follows easily that for a Dirac subspace $L$ one has
$$ker(\omega|_L)= (W_{+}\cap L) \oplus (W_{-}\cap L).$$

The following proposition shows that the integers $k=dim(W_{-}\cap
L)$ and $r=rank(\omega|_L)$ are the only invariants of a Dirac
subspace with respect to the action of the para-Hermitian subgroup
$pH(W)$.
\begin{prop}\label{acttranz} The group $pH(W)$ acts
transitively on the set of Dirac subspaces $L$ with given values
$k,r$.
\end{prop} \noindent{\bf Proof.} Denote $p_{\pm}=F_{\pm}|_L$.
Obviously, $ker\,p_{\pm}=W_{\mp}\cap L$, and we get two linear
spaces
$$L_{\pm}=im\,p_{\pm}\subseteq W_{\pm}\approx L/{W_{\mp}\cap
L}$$ of dimension $n-dim(W_{\mp}\cap L)$. Since $W_{\mp}\cap
L\subseteq ker(\omega|_L)$, we see that the subspaces $L_{\pm}$
have induced, skew-symmetric, bilinear forms $\omega_{\pm}^L$.

Moreover, because of the structure of $ker(\omega|_L)$ as
described above, we see that
$rank\,\omega^{L}_{\pm}=rank\,\omega|_L$ and
$ker\,\omega^L_{\pm}=p_{\pm}(W_{\pm}\cap L)$. (Clearly,
$p_{\pm}|_{W_{\pm}\cap L}$ have kernel zero and are isomorphisms
onto the corresponding images.)

As shown in \cite{C}, it is possible to reconstruct $L$ from each
of the pairs $(L_{\pm},\omega_{\pm}^L)$ namely, $$\begin{array}{l}
L=\{w\in W\;/\;F_{+}(w) \in
L_{+},\;g(F_{-}(w),u)=\omega_{+}^L(F_{+}(w),u),\forall u\in
L_{+}\},\vspace{2mm}\\ L=\{w\in W\;/\;F_{-}(w) \in
L_{-},\;g(F_{+}(w),u)=\omega_{-}^L(F_{-}(w),u),\forall u\in
L_{-}\}.
\end{array}$$

The two formulas have similar justifications, and we refer to the
first only. A straightforward check shows that, for any choice of
a pair $(L_+,\omega_+)$, which consists of an arbitrary subspace
of $W$ and an arbitrary 2-form on that subspace, $L$ defined by
the first formula is an isotropic subspace of $W$, that
$F_{+}(L)=L_{+}$, and that $\omega_{+}$ is induced by $\omega$.
Furthermore, the formula implies
$$ker\,F_{+}|_L= L\cap W_{-}=\{w\in W_{-}\;/\;g(w,u)=0,\;\forall
u\in L_{+}\}, $$ and, since $g|_{W_{-}\times W_{+}}$ is a non
degenerate pairing, $k=n-dim(L_{+})$. Accordingly,
$dim\,L=dim(ker\,F_{+}|_L)+ dim(im\,F_{+}|_L)=n$ and $L$ is the
required Dirac subspace.

Now, if $L,L'$ are Dirac subspaces of $W$ with the same invariants
$k,r$, there exists a transformation $\psi\in Gl(W_{+})$ which
sends the pair $(L_{+},\omega_{+}^L)$ onto
$(L'_{+},\omega_{+}^{L'})$. Obviously, the image of $\psi$ in
$pH(W)$ via the embedding (\ref{groupn}) sends $L$ onto $L'$.
Q.e.d.

Of course, instead of the invariant $k=dim(W_{-}\cap L)$ we may
consider $h=dim(W_{+}\cap L)$. These two numbers are related by
$k+h=n-r$.
%\begin{center}\begin{thebibliography}{xx}
 %\end{center}
\hspace*{7.5cm}{\small \begin{tabular}{l} Department of
Mathematics\\ University of Haifa, Israel\\ E-mail:
vaisman@math.haifa.ac.il \end{tabular}}
\end{document}